\documentclass[a4paper,12pt]{article}
\usepackage{amsmath}
\usepackage{amssymb}
\usepackage{tabularx}
\usepackage{enumerate}
\usepackage{graphicx,subfigure}
\newtheorem{thm}{Theorem}[section]
\newtheorem{lem}[thm]{Lemma}

\newtheorem{defi}{Definition}[section]
\newtheorem{pppp}{Proof}

\newcommand{\qed}{\hspace{1em}\mbox{\raisebox{0.65ex}{\fbox{}}}}

\numberwithin{equation}{section}

\newcommand{\be}{\begin{equation}}
\newcommand{\ee}{\end{equation}}
\newcommand\bes{\begin{eqnarray}} \newcommand\ees{\end{eqnarray}}
\newcommand{\bess}{\begin{eqnarray*}}
\newcommand{\eess}{\end{eqnarray*}}

\newcommand{\bpf}{{\bf Proof:\ \ }}
\newcommand{\epf}{\mbox{}\hfill $\Box$}

\begin{document}

\thispagestyle{empty}

\title{{\bf Analysis of a mutualism
model with stochastic perturbations\thanks{The work is supported by Supported in part by a NSFC Grant No. 11171158,
NSF of Jiangsu Education Committee No. 11KJA110001, "333" Project of Jiangsu Province  Grant No. BRA2011173.}}}
\date{\empty}
\author{Mei Li$^{a, b}$, Hongjun Gao$^{a,}$\footnote{The corresponding author. E-mail address: gaohj@njnu.edu.cn} , Chenfeng Sun$^{b,  c}$, Yuezheng Gong$^{a}$ \\
{\small $^a$ Institute of mathematics, Nanjing Normal University,}\\
{\small Nanjing 210023, PR China}\\
{\small $^b$ School of Applied Mathematics, Nanjing University of Finance and Economics,}\\
{\small Nanjing 210023, PR China}\\
{\small Email: limei@njue.edu.cn }\\
{\small $^c$ Jiangsu Key Laboratory for Numerical Simulation of Large Scale Complex Systems,}\\
{\small Nanjing 210023, PR China}}

\maketitle

\begin{quote}
\noindent
{\bf Abstract.} { 
\small This article is concerned with a mutualism ecological model with stochastic perturbations. The local existence and
uniqueness of a positive solution are obtained with positive initial value, and the asymptotic
behavior to the  problem is studied. Moreover, we show that the solution is stochastically bounded,
uniformly continuous and stochastic permanence. The sufficient conditions for the system to be extinct are given
and the condition for the system to be persistent are also established. At last, some figures are presented to illustrate our main results.
}

\end{quote}

\section{Introduction}

Mutualism is an important biological interaction in nature. It
occurs when one species provides some benefit in exchange for some
benefit, for example, pollinators and flowering plants,  the
pollinators obtain floral nectar (and in some cases pollen) as a
food resource while the plant obtains non-trophic reproductive
benefits through pollen dispersal and seed production. Another
instance is ants and aphids, in which the ants obtain honeydew
food resources excreted by aphids while the aphids obtain
increased survival by the non-trophic service of ant defense
against natural enemies of the aphids. Lots of author have
discussed these models \cite{AR, BJK, CC, CY, GB, HS, HDB, HD1, LT,
TNS}. One of the simplest models is the classical Lotka-Volterra
two-species mutualism model as follow:
\begin{eqnarray}
\left\{
\begin{array}{lll}
\dot{x}(t)=x(t)\big(a_1-b_1x(t)+c_1y(t)\big), \\
\dot{y}(t)=y(t)\big(a_2-b_2y(t)+c_2x(t)\big).
\end{array} \right.
\label{a1}
\end{eqnarray}

Among various types mutualistic model, we should specially mention the following model which was proposed by May (\cite{MRM}) in 1976:
\begin{eqnarray}
\left\{
\begin{array}{ll}
\dot{x}(t)=x(t)\big(r_1-\frac{b_1x(t)}{K_1+y(t)}-\varepsilon_1x(t)\big),& \\
\dot{y}(t)=y(t)\big(r_2-\frac{b_2y(t)}{K_2+x(t)}-\varepsilon_2y(t)\big),&
\end{array} \right.
\label{aode1}
\end{eqnarray}
where $x(t), y(t)$ denote population densities of each species at time t, $r_i, K_i, \alpha_i, \varepsilon_i$ (i=1, 2)
are positive constants, $r_1, r_2$ denote the intrinsic growth rate of species $x(t), y(t)$ respectively, $K_1$ is the capability of species $x(t)$ being short of  $y(t)$, similarly $K_2$ is the capability of species $y(t)$ being short of  $x(t)$. For \eqref{aode1}, there are three trivial equilibrium points
$$ E_1=(0,0),\ \ E_2=(\frac{r_1}{\varepsilon_1+\frac{b_1}{K_1}}, 0),\ \ E_3=(0, \frac{r_2}{\varepsilon_2+\frac{b_2}{K_2}}),$$
and a unique positive interior equilibrium point $E^*=(x^*, y^*)$ satisfies the following equations
\begin{eqnarray}
\left\{
\begin{array}{ll}
r_1-\frac{b_1x(t)}{K_1+y(t)}-\varepsilon_1x(t)=0,& \\
r_2-\frac{b_2y(t)}{K_2+x(t)}-\varepsilon_2y(t)=0,&
\end{array} \right.
\label{aode2}
\end{eqnarray}
where $E^*$ is globally asymptotically stable.

 In addition, population dynamics is inevitably affected by environmental noises,
  May\cite{MR} pointed out the fact that due to environmental fluctuation, the birth rates, carrying capacity,  and
  other parameters involved in the model system exhibit random fluctuation to a greater or lesser extent. Consequently the equilibrium population
   distribution fluctuates randomly around some average values. Therefore lots of authors introduced stochastic perturbation into deterministic models
   to reveal the effect of environmental variability on the population dynamics in mathematical ecology \cite{DS, HWH, JSL, JJS, JJ, LGJM, LLL, LW1, LW2, LW3, MSR, TDHS}.
   So far as our knowledge is concerned, taking into account the effect of randomly fluctuating environment,
   we now add white noise to each equations of the problem \eqref{aode1}.
   Suppose that parameter $r_i$ is stochastically perturbed, with
   $$r_i\rightarrow r_i+\alpha_i\dot{W}_i(t), \  i=1, 2$$
   where $W_1(t), W_2(t)$ are mutually independent Brownian motion, $\alpha_i^2 $ represent the intensities of the white noise. then
   the corresponding deterministic model system \eqref{aode1} may be described by the It$\hat{o}$ problems:
\begin{eqnarray}
\left\{
\begin{array}{ll}
dx(t)=x(t)\big(r_1-\frac{b_1x(t)}{K_1+y(t)}-\varepsilon_1x(t)\big)dt+\alpha_1x(t)d W_1(t),\\
dy(t)=y(t)\big(r_2-\frac{b_2y(t)}{K_2+x(t)}-\varepsilon_2y(t)\big)dt+\alpha_2y(t)d W_2(t).\\
\end{array} \right.
\label{a3}
\end{eqnarray}

In this paper, we will discuss the stability in time average. We now briefly give an outline of the paper.
In the next section, the global existence and uniqueness of the positive solution to problem
\eqref{a3} are proved by using comparison theorem for stochastic equations.
Sections 3 and 4 is devoted to stochastic boundedness, uniformly  H\"{o}lder-continuous.
Section 5 deals with stochastic permanence. Section 6 discusses the persistence in mean and extinction,
sufficient conditions of persistence in mean and extinction are obtained.
Finally in section 7, we carry out numerical simulations to confirm our part results.
throughout this paper, we let $(\Omega, F, {\mathcal\{F_t\}}_{t\geq 0}, P)$ be a complete probability space with a filtration
${\mathcal\{F_t\}}_{t\geq 0}$ satisfying the usual conditions. $X(t)=(x(t), y(t))$ and $|X(t)|=\sqrt{x^2(t)+y^2(t)}.$

We end this section by recalling three definitions and two lemmas which we will use in the forthcoming sections.

\begin{defi} $\cite{LM}$ If for any $0<\varepsilon < 1$, there is a constant $\delta>0$ such that the solution $(x(t), y(t))$ of \eqref{a3} satisfies
$$\limsup_{t\rightarrow\infty}P{|X(t)|}>\delta\}<\varepsilon,$$
for any  initial value $(x_0, y_0)>(0, 0)$, then we say the solution $ X(t)$ be stochastically ultimate boundedness.
\end{defi}

\begin{defi} $\cite{LM}$ If for arbitrary $\varepsilon\in (0,1),$ there are two positive constants $\beta_1$ and $\beta_2$
such that for positive initial data $X_0=(x_0, y_0)$, the solution $X(t)$of problem \eqref{a3} has the property that
$$\liminf_{t\rightarrow \infty}P\{|X(t)|\geq \beta_1\}\geq 1-\varepsilon, \ \   \liminf_{t\rightarrow \infty}P\{|X(t)|\leq \beta_2\}\geq 1-\varepsilon.\ \ $$
Then problem \eqref{a3} is said to be stochastically permanent.
\end{defi}

\begin{defi}$\cite{CLCJ}$ If  $x(t), y(t)$ satisfy the following condition
$$\lim_{t\rightarrow \infty}\frac{1}{t}\int_0^tx(s)ds>0, \ \  \lim_{t\rightarrow \infty}\frac{1}{t}\int_0^tx(s)ds>0 \ \  a.s.$$
The problem of \eqref{a3} is said to be persistence in mean.
\end{defi}
\begin{lem} (Chebyshev's inequality) $\cite{MY}$ If $\delta>0, k>0$ and $X\in L^p(\Omega)$ with $E|X|^k<\infty,$
then,
$$P\{|X|\geq \delta\}\leq \delta^{-k}E|X|^k.$$
\end{lem}

\begin{lem} $\cite{KS}$ Assume that an n-dimensional stochastic process $X(t)$ on  $t\geq 0$
satisfies the condition
$$E|X(t)-X(s)|^\alpha \leq C|t-s|^{1+\beta},\  0\leq s,\, t < +\infty$$
for some positive constants $\alpha, \beta$ and $C$.
There exists a continuous modification $\tilde{X}(t)$ of $X(t)$, which has
the property that for every $\gamma \in (0,\beta/\alpha)$, there is a positive random variable $h(w)$  such that
$$P\big\{\omega: \sup_{0<|t-s|< h(w), 0\leq s, t < +\infty}\frac{|\tilde{X}(t, \omega)-\tilde{X}(s, \omega)|}{|t-s|^\gamma}\leq \frac{2}{1-2^{-\gamma}}\big\}=1.$$
In other words, almost every sample path of $\tilde{X}(t)$ is locally but uniformly H\"{o}lder-continuous with exponent $\gamma$.
\end{lem}

\section{Existence and uniqueness of the positive solution}

First, we show that there exists a unique local  positive solution of  \eqref{a3}.
\begin{lem} For the given positive initial value $(x_0, y_0)$, there is  $\tau \geq 0$ such that
 problem \eqref{a3}
admits a unique positive local solution $(x(t), y(t))$  a.s. for $t\in [0, \tau).$
\end{lem}
\bpf
  We first set a change of variables : $u(t)=\ln x(t), v(t)=\ln y(t)$, then problem \eqref{a3} deduces to
\begin{eqnarray}
\left\{
\begin{array}{lll}
du(t)=\big(r_1-\alpha_1^2/2-\frac{b_1e^{u(t)}}{K_1+e^{v(t)}}-\varepsilon_1 e^{u(t)}\big)dt+\alpha_1 dW_1(t),\\
dv(t)=\big(r_2-\alpha_2^2/2-\frac{b_2e^{v(t)}}{K_2+e^{u(t)}}-\varepsilon_2 e^{v(t)}\big)dt+\alpha_2 dW_2(t)\\
\end{array} \right.
\label{Hb}
\end{eqnarray}
on $t\geq 0$ with initial value $u(0)=\ln{x_0}, v(0)=\ln{y_0} $.
Obviously, the coefficients of \eqref{Hb} satisfy the local Lipschitz condition,
then, making use of the theorem \cite{FA, MX} about existence and uniqueness for stochastic differential equation
there is a unique local solution $(u(t), v(t))$ on $t \in [0,\tau )$, where $\tau$ is the explosion time. Hence,
by It$\hat{o}$'s formula, $(x(t), y(t))$ is a unique positive local solution to problem \eqref{a3} with positive initial value.

Next we need to prove solution is  global, that is $\tau=\infty$.
\begin{thm} For any positive initial value $(x_0, y_0)$, there exists a unique global positive solution $(x(t), y(t))$ to problem \eqref{a3}, which
satisfies
$$\lambda (t)\leq x(t)\leq \Lambda(t), \  \theta (t)\leq y(t)\leq \Theta(t), \ t\geq 0, \  a.s.$$
where $\lambda,\  \Lambda,\  \theta$ and $\Theta$ are defined as \eqref{HC}, \eqref{12}, \eqref{13} and \eqref{14}.
\end{thm}
\bpf   \cite{JJS}was the main source of inspiration for its proof. Because of $(x(t), y(t))$ is positive, from the first equation of \eqref{a3} we get
$$dx(t)\leq x(t)\big(r_1-\varepsilon_1x(t)\big)dt+\alpha_1x(t)d W_1(t).$$
Define the following problem
\begin{eqnarray}
\left\{
\begin{array}{l}
d\Lambda(t)=\Lambda(t)\big(r_1-\varepsilon_1\Lambda(t)\big)dt+\alpha_1\Lambda(t) dW_1(t),\\
\Lambda(0)=x_0,
\end{array} \right.
\label{11}
\end{eqnarray}
then
\begin{eqnarray}
\Lambda(t)=\frac{\exp\big((r_1-\frac{\alpha_1^2}{2})t+\alpha_1W_1(t)\big)}{\frac{1}{x_0}+\varepsilon_1\int^t_0\exp\big((r_1-\frac{\alpha_1^2}{2})s+\alpha_1W_1(s)\big)ds}
\label{HC}
\end{eqnarray}
is the unique solution of \eqref{11}, and it follows from the comparison theorem for stochastic equations that
$$x(t)\leq \Lambda(t),  \  t\in [0,\tau), \ a.s.$$
On the other hand,
\begin{eqnarray*}
\begin{array}{lll}
dx(t)&\geq &x(t)\big(r_1-\frac{b_1x(t)}{K_1}-\varepsilon_1x(t)\big)dt+\alpha_1dW_1(t) \\
&=&x(t)\big(r_1-(\frac{b_1}{K_1}+\varepsilon_1)x(t)\big)dt+\alpha_1dW_1(t).
\end{array}
\end{eqnarray*}
Obviously,
\begin{eqnarray}
\lambda(t)=\frac{\exp((r_1-\frac{\alpha_1^2}{2})t+\alpha_1W_1(t))}{\frac{1}{x_0}+(\frac{b_1}{K_1}+\varepsilon_1)
\int^t_0\exp((r_1-\frac{\alpha_1^2}{2})s+\alpha_1W_1(s))ds}
\label{12}
\end{eqnarray}
is the solution to the problem
\begin{eqnarray}
\left\{
\begin{array}{llll}
d\lambda(t)=\lambda(t)\big(r_1-(\frac{b_1}{K_1}+\varepsilon_1)\lambda(t)\big)dt+\alpha_1\lambda(t) dW_1(t),\\
\lambda(0)=x_0,
\end{array} \right.
\label{Hd}
\end{eqnarray}
and
$$x(t)\geq \lambda(t),  \  t\in [0,\tau), \ a.s.$$
Similarly, we can get
$$y(t)\leq \Theta(t),  \  t\in [0,\tau), \  a.s,$$
where
\begin{eqnarray}
\Theta(t)=\frac{\exp\big((r_2-\frac{\alpha_2^2}{2})t+\alpha_2W_2(t)\big)}{\frac{1}{y_0}+\varepsilon_2\int^t_0
\exp((r_2-\frac{\alpha_2^2}{2})s+\alpha_2W_2(s))ds},
\label{13}
\end{eqnarray}
and,
$$y(t)\geq \theta(t),  \  t\in [0,\tau), \ a.s.$$
where
\begin{eqnarray}
\theta(t)=\frac{\exp((r_2-\frac{\alpha_2^2}{2})t+\alpha_2W_2(t))}{\frac{1}{y_0}+(\frac{b_2}{K_2}+\varepsilon_2)\int^t_0
\exp((r_2-\frac{\alpha_2^2}{2})s+\alpha_2W_2(s))ds}.
\label{14}
\end{eqnarray}
Combining \eqref{13} and \eqref{14},  we obtain
$$\lambda (t)\leq x(t)\leq \Lambda(t), \  \theta (t)\leq y(t)\leq \Theta(t), t\geq 0, \  a.s.$$
Since that $\Lambda(t), \lambda (t), \theta (t)$ and $\Theta(t)$ exist for any $t>0$, it follows from the comparison theorem for stochastic equations \cite{IW}
that $(x(t), y(t))$ exists globally.
\epf

\section{Stochastically ultimate boundedness }
  In a population dynamical system, the nonexplosion property is often not good enough but the property of ultimate boundedness is more desired.
Now, let us present a theorem  about the Stochastically ultimate boundedness of \eqref{a3} for any positive initial value.

\begin{thm} For any positive initial value $(x_0, y_0)$, the solution $X(t)$ of problem \eqref{a3} is stochastically ultimate boundedness.
\end{thm}\
\bpf As in \cite{LM} we define the function $U=e^tx^k,\ k> 0$.  By the It$\hat{o}$'s formula:
\begin{eqnarray*}
\begin{array}{lll}
d(e^tx^k)&=&e^tx^kdt+ke^tx^{k-1}dx+\frac{1}{2}k(k-1)e^tx^{k-2}(dx)^2\\
&=&e^tx^kdt+ke^tx^{k-1}x(r_1-\frac{b_1x}{K_1+y}-\varepsilon_1x)dt+\alpha_1ke^tx^kdW_1(t)\\
& &+\frac{1}{2}k(k-1)e^tx^{k-1}\alpha_1^2x^2dt\\
&=&e^tx^k[1+k(r_1-\frac{b_1x}{K_1+y}-\varepsilon_1x)+\frac{1}{2}\alpha_1^2k(k-1)]dt+\alpha_1ke^tx^kdW_1(t).
\end{array}
\end{eqnarray*}
Application of  young's inequality yields,
\begin{eqnarray*}
\begin{array}{ll}
&e^tx^k[1+k(r_1-\frac{b_1x}{K_1+y}-\varepsilon_1x)+\frac{1}{2}\alpha_1^2k(k-1)]\\
&\leq e^t[(1+kr_1+\frac{1}{2}\alpha_1^2k(k-1))x^k-k\varepsilon_1x^{k+1}]\\
&\leq \frac{e^t}{\varepsilon_1^k}\big[\frac{1+kr_1+\frac{k(k-1)}{2}\alpha_1^2}{k+1}\big]^{k+1}:=H_1(k)e^t,
\end{array}
\end{eqnarray*}
where  $H_1(k)=\frac{1}{\varepsilon_1^k}\big[\frac{1+kr_1+\frac{k(k-1)}{2}\alpha_1^2}{k+1}\big]^{k+1}$. Therefore,
$$d(e^tx^k)\leq H_1(k)e^tdt+\alpha_1ke^tx^kdW_1(t),$$
Taking expectation to obtain
$$E(e^tx^k)-E(x_0^k)\leq H_1(k)e^t.$$
Thus,
\begin{eqnarray}
 \limsup_{t\rightarrow\infty} Ex^k\leq H_1(k).
 \label{31}
 \end{eqnarray}
Similarly, we have
\begin{eqnarray}
\limsup_{t\rightarrow\infty}Ey^k\leq H_2(k),
\label{32}
\end{eqnarray}
where $H_2(k)=\frac{1}{\varepsilon_2^k}\big[\frac{1+kr_2+\frac{k(k-1)}{2}\alpha_2^2}{k+1}\big]^{k+1}.$
We now combine \eqref{31}, \eqref{32} and  the formula  $\big[x(t)^2+y(t)^2\big]^{\frac{k}{2}}\leq 2^{\frac{k}{2}}\big[x(t)^k+y(t)^k\big]$ to yield
$$\limsup_{t\rightarrow\infty} E|X|^k\leq 2^{\frac{k}{2}}(H_1(k)+H_2(k))<+\infty.$$
By the Lemma 1.1 we can complete the proof.
\epf

\section{Uniformly  H\"{o}lder-continuous }
Now, let us discuss the uniformly  H\"{o}lder-continuous about the positive solution of  problem \eqref{a3}.
\begin{thm} Let $X(t)$ be a positive solution of  problem \eqref{a3} for any positive initial value $X(0)=(x_0, y_0)$, almost every sample path
of $X(t)$ to  \eqref{a3} is uniformly H\"{o}lder-continuous.
\end{thm}
\bpf The proof is motivated by the arguments in \cite{LW2}.
The first equation of \eqref{a3} is equivalent to the  following stochastic integral equation
$$x(t)=x_0+\int_0^tx(s)[r_1-\frac{b_1x(s)}{K_1+y(s)}-\varepsilon_1x(s)]ds+\int_0^t\alpha_1x(s)dW_1(s).$$
By Theorem 3.2 and the inequality \eqref{31},
we have
\begin{eqnarray*}
\begin{array}{ll}
&E\big|x(s)[r_1-\frac{b_1x(s)}{K_1+y(s)}-\varepsilon_1x(s)]\big|^k \\
&\leq 0.5 Ex^{2k}(s)+0.5E(r_1-\frac{b_1x(s)}{K_1+y(s)}-\varepsilon_1x(s))^{2k}\\
&\leq 0.5Ex^{2k}(s)+2^{2k-2}r_1^{2k}+2^{2k-2}Ex^{2k}(s)\\
&\leq 0.5H_1{(2k)}+2^{2k-2}r_1^{2k}+2^{2k-2}H_1{(2k)}\\
&=: H_{11}(k).
\end{array}
\end{eqnarray*}
and
$$E|\alpha_1x(s)|^k\leq \alpha_1^kH_1(k).$$
Using the moment inequality (\cite{MX}) gives that
\begin{eqnarray*}
\begin{array}{lll}
&E|\int_{t_1}^{t_2}\alpha_1x(s)dW_1(s)|^k\\
&\leq \alpha_1^k[0.5k(k-1)]^{0.5k}(t_2-t_1)^{0.5(k-2)}\int_{t_1}^{t_2}Ex^k(s)ds\\
&\leq \alpha_1^k[0.5k(k-1)]^{0.5k}(t_2-t_1)^{0.5(k)}H_1(k)
\end{array}
\end{eqnarray*}
for $0\leq t_1\leq t_2$ and $k>2.$
Therefore we obtain
\begin{eqnarray*}
\begin{array}{lll}
& &E|x(t_2)-x(t_1)|^k\\
&=&E\big|\int_{t_1}^{t_2}x(s)[r_1-\frac{b_1x(s)}{K_1+y(s)}-\varepsilon_1x(s)]ds+\int_{t_1}^{t_2}\alpha_1x(s)dW_1(s)\big|^k\\
&\leq& 2^{k-1}E\big|\int_{t_1}^{t_2}x(s)[r_1-\frac{b_1x(s)}{K_1+y(s)}-\varepsilon_1x(s)]ds\big|^k\\
& &+2^{k-1}E\big|\int_{t_1}^{t_2}\alpha_1x(s)dW_1(s)\big|^k\\
&\leq &2^{k-1}(\int_{t_1}^{t_2}ds)^{k-1}\int_{t_1}^{t_2}E|x(s)[r_1-\frac{b_1x(s)}{K_1+y(s)}-\varepsilon_1x(s)]|^kds\\
& &+2^{k-1}\alpha_1^k[0.5k(k-1)(t_2-t_1)]^{0.5k}H_1(k)\\
& \leq &2^{k-1}(t_2-t_1)^kH_{11}(k)+2^{k-1}\alpha_1^k[0.5k(k-1)(t_2-t_1)]^{0.5k}H_1(k)\\
&\leq& 2^{k-1}(t_2-t_1)^{0.5k}\big\{(t_2-t_1)^{0.5k}+[0.5k(k-1)]^{0.5k}\big\}H(k),
\end{array}
\end{eqnarray*}
for $0<t_1<t_2<\infty, t_2-t_1\leq 1, k>2,$ where $H(k)=\max\big\{H_{11}(k), \alpha_1^kH_1(k) \big\}.$  Hence, it follows from Lemma 1.2 that
almost every sample path of $x(t)$ is locally but uniformly H\"{o}lder continuous with exponent $\gamma\in (0, \frac{k-2}{2k}),$ and almost
every sample path of $x(t)$ is uniformly continuous on $t\geq 0$.
Similarly, almost every sample path of $y(t)$ is uniformly continuous on $t\geq 0$. All in all, almost every sample path of $X(t)=(x(t), y(t))$
of \eqref{a3} be  uniformly continuous on $t\geq 0$.
\epf

\section{Stochastic permanence}

In the study of population models, permanence is one of the most interesting and important topics. We will discuss the property
by using the method as in \cite{LW3} in this section.

\begin{thm} Problem \eqref{a3} is stochastically permanent.
\end{thm}
\bpf For a positive constant $\eta <1$, we set a function
$$U(X)=\frac{1}{\eta}(1+\frac{1}{x})^\eta + \frac{1}{\eta}(1+\frac{1}{y})^\eta, $$
Straightforward compute $d U(X)$ by  It$\hat{o}$ formula
\begin{eqnarray*}
\begin{array}{lll}
& &d U(X)\\
&=&(1+\frac{1}{x})^{\eta-1}d(\frac{1}{x})+(1+\frac{1}{y})^{\eta-1}d(\frac{1}{y})\\
& &+0.5(\eta-1)(1+\frac{1}{x})^{\eta-2}(d(\frac{1}{x}))^2+0.5(\eta-1)(1+\frac{1}{y})^{\eta-2}(d(\frac{1}{y}))^2\\
&=&(1+\frac{1}{x})^{\eta-2}\big\{ (1+\frac{1}{x})[-\frac{1}{x}(r_1-\frac{b_1x(t)}{K_1+y(t)}-\varepsilon_1x(t))+0.5(\eta-1)\alpha_1^2\frac{1}{x^2}]\big\}dt\\
& &+(1+\frac{1}{x})^{\eta-2}\big\{ (1+\frac{1}{x})[-\frac{1}{x}(r_1-\frac{b_1x(t)}{K_1+y(t)}-\varepsilon_1x(t))+0.5(\eta-1)\alpha_1^2\frac{1}{x^2}]\big\}dt\\
& &-\frac{1}{x}(1+\frac{1}{x})^{\eta-1}\alpha_1dW_1(t)-\frac{1}{y}(1+\frac{1}{y})^{\eta-1}\alpha_2dW_2(t)\\
&\leq& (1+\frac{1}{x})^{\eta-2}\big\{-\frac{1}{x^2}[r_1-0.5(\eta-1)\alpha_1^2]+\frac{1}{x}(-r_1+\varepsilon_1+
\frac{b_1}{K_1})+(\varepsilon_1+\frac{b_1}{K_1})\big\}dt\\
& &+(1+\frac{1}{y})^{\eta-2}\big\{-\frac{1}{y^2}[r_1-0.5(\eta-1)\alpha_1^2]+\frac{1}{y}(-r_2+\varepsilon_2+
\frac{b_2}{K_2})+(\varepsilon_2+\frac{b_2}{K_2})\big\}dt\\
& &-\frac{1}{x}(1+\frac{1}{x})^{\eta-1}\alpha_1dW_1(t)-\frac{1}{y}(1+\frac{1}{y})^{\eta-1}\alpha_2dW_2(t).
\end{array}
\end{eqnarray*}
Let us choose $\mu$ sufficiently small to satisty
$$0< \frac{\mu}{\eta} <  \min\{r_1, r_2\}$$
Set $V(X)=e^{\mu t}U(X),$ applying the It$\hat{o}^,s$ formula, we obtain
\begin{eqnarray*}
\begin{array}{lll}
& &d V(X)=\mu e^{\mu t}U(X)dt+e^{\mu t}d U(X)\\
&\leq& \mu e^{\mu t}[\frac{1}{\eta}(1+\frac{1}{x})^\eta + \frac{1}{\eta}(1+\frac{1}{y})^\eta]dt+e^{\mu t}d U(X)\\
&\leq& e^{\mu t}(1+\frac{1}{x})^{\eta-2}\{-\frac{1}{x^2}[r_1-0.5(\eta-1)\alpha_1^2-\frac{\mu}{\eta}]+\frac{1}{x}(-r_1+\frac{2\mu}{\eta})\\
& &+(1+\frac{1}{x})(\varepsilon_1+\frac{b_1}{K_1})+\frac{\mu}{\eta}\}dt\\
& &+e^{\mu t}(1+\frac{1}{y})^{\eta-2}\{-\frac{1}{y^2}[r_2-0.5(\eta-1)\alpha_2^2-\frac{\mu}{\eta}]+\frac{1}{y}(-r_2+\frac{2\mu}{\eta})\\
& &+(1+\frac{1}{y})(\varepsilon_2+\frac{b_2}{K_2})+\frac{\mu}{\eta}\}dt \\
& &-e^{\mu t}[\frac{1}{x}(1+\frac{1}{x})^{\eta-1}\alpha_1dW_1(t)+\frac{1}{y}(1+\frac{1}{y})^{\eta-1}\alpha_2dW_2(t)]\\
&\leq& e^{\mu t}[(1+\frac{1}{x})^{\eta-1}(\varepsilon_1+\frac{b_1}{K_1}+\frac{\mu}{\eta})+(1+\frac{1}{y})^{\eta-1}(\varepsilon_2+\frac{b_2}{K_2}+\frac{\mu}{\eta})]dt\\
& &-e^{\mu t}[\frac{1}{x}(1+\frac{1}{x})^{\eta-1}\alpha_1dW_1(t)+\frac{1}{y}(1+\frac{1}{y})^{\eta-1}\alpha_2dW_2(t)]\\
&\leq& L_1e^{\mu t}dt-e^{\mu t}[\frac{1}{x}(1+\frac{1}{x})^{\eta-1}\alpha_1dW_1(t)+\frac{1}{y}(1+\frac{1}{y})^{\eta-1}\alpha_2dW_2(t)],
\end{array}
\end{eqnarray*}
where $L_1=\max \{\varepsilon_1+\frac{b_1}{K_1}+\frac{\mu}{\eta},  \varepsilon_2+\frac{b_2}{K_2}+\frac{\mu}{\eta} \}.$

Integrating and then taking expectations yields
$$E[V(X)]=e^{\mu t}E(U(X))\leq \frac{1}{\eta}(1+\frac{1}{x_0})^\eta + \frac{1}{\eta}(1+\frac{1}{y_0})^\eta+\frac{L_1}{\mu}(e^{\mu t}-1).$$
Therefore,
$$\limsup_{t\rightarrow+\infty}E[\frac{1}{x^\eta(t)}]\leq \limsup_{t\rightarrow+\infty}E[(1+\frac{1}{x(t)})^\eta+(1+\frac{1}{y(t)})^\eta]\leq \frac{\eta L_1}{\mu}=L,$$
and
$$\limsup_{t\rightarrow+\infty}E[\frac{1}{y^\eta(t)}]\leq \limsup_{t\rightarrow+\infty}E[(1+\frac{1}{x(t)})^\eta+(1+\frac{1}{y(t)})^\eta]\leq \frac{\eta L_1}{\mu}=L.$$
For arbitrary $\varepsilon\in(0,1)$, choose $\beta_1(\varepsilon)=(\frac{\varepsilon}{L})^{\frac{1}{\eta}},$  we yield the following inequality making Chebyshev's inequality,
$$P\{x(t)<\beta_1\}=P\{\frac{1}{x^\eta(t)}>\frac{1}{\beta_1^\eta}\}\leq \frac{E[\frac{1}{x^\eta(t)}]}{\beta_1^{-\eta}},$$
$$P\{y(t)<\beta_1\}=P\{\frac{1}{y^\eta(t)}>\frac{1}{\beta_1^\eta}\}\leq \frac{E[\frac{1}{y^\eta(t)}]}{\beta_1^{-\eta}}.$$
Hence,
$$\limsup_{t\rightarrow+\infty}P\{|X(t)|<\beta_1\}\leq \beta_1^\eta L=\varepsilon,$$
Then,
$$\liminf_{t\rightarrow+\infty}P\{|X(t)|\geq\beta_1\}\geq 1-\varepsilon.$$
Using Chebyshev's inequality and \eqref{31}, \eqref{32} we can prove that\\
for arbitrary $\varepsilon\in(0,1)$, there is a positive constant $\beta_2$ such that
$$\liminf_{t\rightarrow+\infty}P\{|X(t)|\leq \beta_2\}\geq 1-\varepsilon.$$
\epf

\section{Persistence in mean and extinction}
In the description of population dynamics, it is critical to discuss the property of persistence in mean and extinction.

\begin{thm} Suppose that $r_i>\frac{\alpha_i^2}{2},  (i=1, 2), $ $X(t)$ is the positive solution to \eqref{a3} with positive initial value
$(x_0, y_0)$, then the problem \eqref{a3} is persistent in mean.
\end{thm}
\bpf The method is similar to \cite{JJS}. We first deduce
$$\lim_{t \rightarrow \infty} \frac{\ln x(t)}{t}=0, \ \lim_{t\rightarrow \infty}\frac{\ln y(t)}{t}=0.$$
For $t\geq T >0 $, we compute from \eqref{HC} that
$$\Lambda(t)=\frac{\exp\big((r_1-\frac{\alpha_1^2}{2})(t-T)+\alpha_1(W_1(t)-W_1(T))\big)}{\frac{1}{x(T)}+\varepsilon_1\int^t_T
\exp\big((r_1-\frac{\alpha_1^2}{2})s+\alpha_1W_1(s)\big)ds}.$$
Hence,
\begin{eqnarray*}
\begin{array}{lll}
& &\frac{1}{\Lambda(t)}
=\frac{\frac{1}{x(T)}+\varepsilon_1\int^t_T\exp\big((r_1-\frac{\alpha_1^2}{2})s+\alpha_1W_1(s)\big)ds}{\exp\big((r_1-\frac{\alpha_1^2}{2})(t-T)+
\alpha_1(W_1(t)-W_1(T))\big)}\\
& \geq& \exp\big( -(r_1-\frac{1}{2}\alpha_1^2)(t-T)-\alpha_1(W_1(t)-W_1(T))\big)\times\\
& &\big[\frac{1}{x(T)}+\varepsilon_1\int^t_T \exp((r_1-\frac{\alpha_1^2}{2})v+\alpha_1W_1(v))dv \big]\\
&\geq & \frac{\varepsilon_1}{r_1-\frac{\alpha_1^2}{2}}\exp\big( (r_1-\frac{\alpha_1^2}{2})T+\alpha_1W_1(T)\big)\big(1-\exp(- (r_1-\frac{\alpha_1^2}{2})(t-T))\big)\times\\
& &\exp\big(\alpha_1(\min_{0\leq v \leq t}W_1(v)-\max_{0\leq v \leq t}W_1(v))\big)\\
&=:&H_1(t)\exp\big(\alpha_1(\min_{0\leq v \leq t}W_1(v)-\max_{0\leq v \leq t}W_1(v))\big),
\end{array}
\end{eqnarray*}
where $H_1(t)=\frac{\varepsilon_1}{r_1-\frac{\alpha_1^2}{2}}e^{(r_1-\frac{\alpha_1^2}{2})T+\alpha_1W_1(T)}\big(1-e^{- (r_1-\frac{\alpha_1^2}{2})(t-T)}\big ).$\\
Taking logarithm to obtain
$$-\ln\Lambda(t)\geq \ln H_1(t)+ \alpha_1\big(\min_{0\leq v \leq t}W_1(v)-\max_{0\leq v \leq t}W_1(v)\big),$$
which implies that
$$\frac{\ln\Lambda(t)}{t}\leq -\frac{\ln H_1(t)}{t}-\alpha_1\frac{\min_{0\leq v \leq t}W_1(v)-\max_{0\leq v \leq t}W_1(v)}{t}.$$
The distributions of $\max_{0\leq v \leq t}W_1(v),$ is same as $|W_1(t)|$, $\min_{0\leq v \leq t}W_1(v)$ have same distribution
as $-\max_{0\leq v \leq t}W_1(v).$
Moreover, $\frac{\ln H_1(t)}{t}\rightarrow \infty$ as $t\rightarrow \infty,$
by the strong law of large numbers we get $\lim_{t\rightarrow \infty}\sup\frac{\ln\Lambda(t)}{t}\leq 0.$ Then
$$\lim_{t\rightarrow \infty}\sup\frac{\ln x(t)}{t}\leq 0.$$
On the other hand, from the \eqref{12} we have:
$$\frac{1}{\lambda(t)}=\frac{1}{x_0}e^{-(r_1-\frac{\alpha_1^2}{2})t-\alpha_1W_1(t)}+(\frac{b_1}{K_1}+\varepsilon_1)
\int_0^te^{-(r_1-\frac{\alpha_1^2}{2})(t-s)-\alpha_1(W_1(t)-W_1(s))}ds$$
$$\leq e^{\alpha_1(\max_{0\leq s \leq t}W_1(s)-W_1(t))}\big[\frac{1}{x_0}e^{-(r_1-\frac{\alpha_1^2}{2})t}+(\frac{b_1}{K_1}+\varepsilon_1)\int_0^te^{-(r_1-\frac{\alpha_1^2}{2})(t-s)}ds\big].$$
Similarly, we can deduce
$$\frac{1}{\lambda(t)}\geq e^{\alpha_1(\min_{0\leq s \leq t}W_1(s)-W_1(t))}\big[\frac{1}{x_0}e^{-(r_1-\frac{\alpha_1^2}{2})t}+(\frac{b_1}{K_1}+\varepsilon_1)\int_0^te^{-(r_1-\frac{\alpha_1^2}{2})(t-s)}ds\big].$$
Note that $$\eta(t)=\frac{1}{\frac{1}{x_0}e^{-(r_1-\frac{\alpha_1^2}{2})t}+(\frac{b_1}{K_1}+\varepsilon_1)\int_0^te^{-(r_1-\frac{\alpha_1^2}{2})(t-s)}ds}$$
is the solution of the problem
\begin{eqnarray}
\left\{
\begin{array}{llll}
\dot{\eta}(t)=\eta(t)\big(r_1-\frac{\alpha_1^2}{2}-(\frac{b_1}{K_1}+\varepsilon_1)\eta(t)\big),\\
\eta(0)=x_0,
\end{array} \right.
\end{eqnarray}
we have
$$e^{\alpha_1(\min_{0\leq s \leq t}W_1(s)-W_1(t))}\frac{1}{\eta(t)}\leq \frac{1}{\lambda(t)} \leq e^{\alpha_1(\max_{0\leq s \leq t}W_1(s)-W_1(t))}\frac{1}{\eta(t)},$$
that is
$$\alpha_1(W_1(t)-\max_{0\leq s \leq t}W_1(s))\leq \ln\lambda(t)-\ln\eta(t)\leq \alpha_1(W_1(t)-\min_{0\leq s \leq t}W_1(s)).$$
Making use of the large number theorem and the distribution of $|W_1(t)|$, we get
$$\lim_{t\rightarrow \infty}\frac{\ln \lambda(t)}{t}= 0,$$
Therefore
$$\lim_{t\rightarrow \infty}\sup\frac{\ln x(t)}{t}\geq 0.$$
Hence,
$$\lim_{t\rightarrow \infty}\frac{\ln x(t)}{t}=0.$$
Similarly, we yield that
$$\lim_{t\rightarrow \infty}\frac{\ln y(t)}{t}=0.$$
Integrating the first equation of (2.1) from $0$ to $t$, we yield
$$b_1\int_0^t\frac{x(s)}{K_1+y(s)}ds=-(\ln x(t)-\ln x_0)+(r_1-\frac{\alpha_1^2}{2})t+\alpha_1W_1(t)-\varepsilon_1\int_0^tx(s)ds,$$
because of       $\int_0^tx(s)ds \geq K_1\int_0^t\frac{x(s)}{K_1+y(s)}ds,$ we obtain

$$b_1\frac{1}{t}\int_0^tx(s)ds \geq K_1\frac{1}{t}\big[-(\ln x(t)-\ln x_0)+(r_1-\frac{\alpha_1^2}{2})t+\alpha_1W_1(t)-\varepsilon_1\int_0^tx(s)ds\big],$$
that is
$$(b_1+\varepsilon_1K_1)\frac{1}{t}\int_0^tx(s)ds \geq -K_1\frac{\ln x(t)-\ln x_0}{t}+K_1(r_1-\frac{\alpha_1^2}{2})+K_1\alpha_1\frac{W_1(t)}{t}.$$
Since that $\lim_{t\rightarrow\infty}\frac{W_1(t)}{t}=0$, and $\lim_{t \rightarrow \infty} \frac{\ln x(t)}{t}=0$, we get
$$\lim_{t \rightarrow \infty} \frac{\int_0^tx(s)ds}{t}\geq \frac{K_1(r_1-\frac{\alpha_1^2}{2})}{b_1+\varepsilon_1K_1}> 0, \ a. s. $$
Similarly, we yield
$$\lim_{t \rightarrow \infty} \frac{\int_0^ty(s)ds}{t}\geq \frac{K_2(r_2-\frac{\alpha_2^2}{2})}{b_2+\varepsilon_2K_2}> 0, \ a. s. $$
The proof is completed
\epf
\begin{thm} Let $X(t)=(x(t), y(t))$ be a positive solution of  \eqref{a3} with positive initial value $X(0)=(x_0, y_0)$, then\\
$(A)$ If $r_1<\frac{\alpha_1^2}{2}, r_2>\frac{\alpha_2^2}{2}$, then $x(t)$ is extinction, $y(t)$ is persistent in mean.\\
$(B)$ If $r_1>\frac{\alpha_1^2}{2}, r_2<\frac{\alpha_2^2}{2}$, then $y(t)$ is extinction, $x(t)$ is persistent in mean.\\
$(C)$ If $r_1<\frac{\alpha_1^2}{2}, r_2<\frac{\alpha_2^2}{2}$, then $x(t), y(t)$ be extinction.
 \end{thm}
\bpf We first prove part $(A)$ of the theorem. The proof of $(B), (C)$ is similar.
It follows from the first equation of (2.1) that
$$du(t)\leq (r_1-\frac{\alpha_1^2}{2})dt+\alpha_1dW_1(t).$$
Apply the comparison theorem for stochastic equations and the diffusion processes, we deduce that
$$\lim_{t\rightarrow\infty}u(t)=-\infty,$$
i.e.\\
$$\lim_{t\rightarrow\infty}x(t)=0,  \   a.s.$$
Hence for any small $\varepsilon>0$, there exist $t_0$ and a set $\Omega_\varepsilon$ such that $P(\Omega_\varepsilon)\geq 1-\varepsilon$ and
$\frac{x(t)}{K_2+x(t)}\leq \varepsilon$ for $t \geq t_0$ and $\omega \in \Omega_\varepsilon.$
Therefore, the second equation of \eqref{a3} becomes
$$dy(t)=y(t)\big(r_2-\frac{b_2y(t)}{K_2+x(t)}-\varepsilon_2y(t)\big)dt+\alpha_2y(t)d W_2(t)$$
$$=y(t)\big(r_2-\frac{b_2y(t)}{K_2}+\frac{b_2y(t)x(t)}{K_2(K_2+x(t))}-\varepsilon_2y(t)\big)dt+\alpha_2y(t)d W_2(t).$$
We can yield
$$y(t)\geq y(t)\big(r_2-(\frac{b_2}{K_2}+\varepsilon_2)y(t)\big)dt+\alpha_2y(t)d W_2(t),$$
$$y(t) \leq y(t)\big(r_2-(\frac{b_2}{K_2}(1-\varepsilon)+\varepsilon_2)y(t)\big)dt+\alpha_2y(t)d W_2(t).$$
If $r_2>\frac{\alpha_2^2}{2},$ using comparison theorem for stochastic equations, we get
$$\lim_{t \rightarrow \infty}\inf \frac{\int_0^ty(s)ds}{t}\geq \frac{r_2-\frac{\alpha_2^2}{2}}{\frac{b_2}{K_2}+\varepsilon_2},
 \  \lim_{t \rightarrow \infty} \sup\frac{\int_0^ty(s)ds}{t}\leq\frac{r_2-\frac{\alpha_2^2}{2}}{\frac{b_2(1-\varepsilon)}{K_2}+\varepsilon_2},$$
which implies that
$$\lim_{t \rightarrow \infty}\frac{\int_0^ty(s)ds}{t}= \frac{r_2-\frac{\alpha_2^2}{2}}{\frac{b_2}{K_2}+\varepsilon_2}>0, \ \ a.s.$$
That is  $y(t)$ is persistent in mean.\\

\epf
\section{Numerical simulations}

Now let us make use of Milstein's  method \cite{HD} to illustrate the analytical findings. Consider the following discretization system:
\begin{eqnarray}
\left\{
\begin{array}{lll}
x^{(n+1)}&=&x^{(n)}+x^{(n)}\big(r_1-\frac{b_1x^{(n)}}{K_1+y^{(n)}}-\varepsilon_1 x^{(n)}\big)\Delta t+\alpha_1 x^{(n)}\sqrt{\Delta t}\xi^{(n)}\\
& &  +\frac{\alpha_1^2}{2} x^{(n)}((\xi^{(n)})^2-1)\Delta t,\\
y^{(n+1)}&=&y^{(n)}+y^{(n)}\big(r_2-\frac{b_2y^{(n)}}{K_2+x^{(n)}}-\varepsilon_2 y^{(n)}\big)\Delta t+\alpha_2 y^{(n)}\sqrt{\Delta t}\zeta^{(n)}\\
  & &      +\frac{\alpha_2^2}{2} y^{(n)}((\zeta^{(n)})^2-1)\Delta t,
\end{array} \right.
\label{d3}
\end{eqnarray}
where $\xi^{(n)}$ and $\zeta^{(n)},$ n=1, 2, ..., N, are the
Gaussian random variables $N(0, 1)$. In Figure 1 we choose
$r_1=1.2, r_2=1, \varepsilon_1=0.8, \varepsilon_2=0.7, b_1=0.7,
b_2=0.9, K_1=K_2=2$, step size $\Delta$ t=0.001. The only
difference between conditions of Figure 1 is that the value of
$\alpha_1, \alpha_2$ . In Figure 1 (a), we choose $\alpha_1=
\alpha_2=0$, we can see the positive equilibrium point $E^*$ is
globally stable;
 In Figure 1 (b)-(d)we choose the value of $\alpha_1, \alpha_2$ such that $r_1<\frac{\alpha_1^2}{2}, r_2<\frac{\alpha_2^2}{2};
  r_1>\frac{\alpha_1^2}{2}, r_2<\frac{\alpha_2^2}{2}; r_1>\frac{\alpha_1^2}{2}, r_2>\frac{\alpha_2^2}{2}$ respectively,
then in view Theorem 6.2, we confirms them.  By comparing In
Figure 1 (a), (b), (c), we can see that small random perturbation
can retain the stochastic system permanent; sufficiently large
random perturbation leads to the stochastic system extinct.
\begin{figure}[!hbp]
\centering
 (a)\subfigure{
      \includegraphics[width=6.5cm]{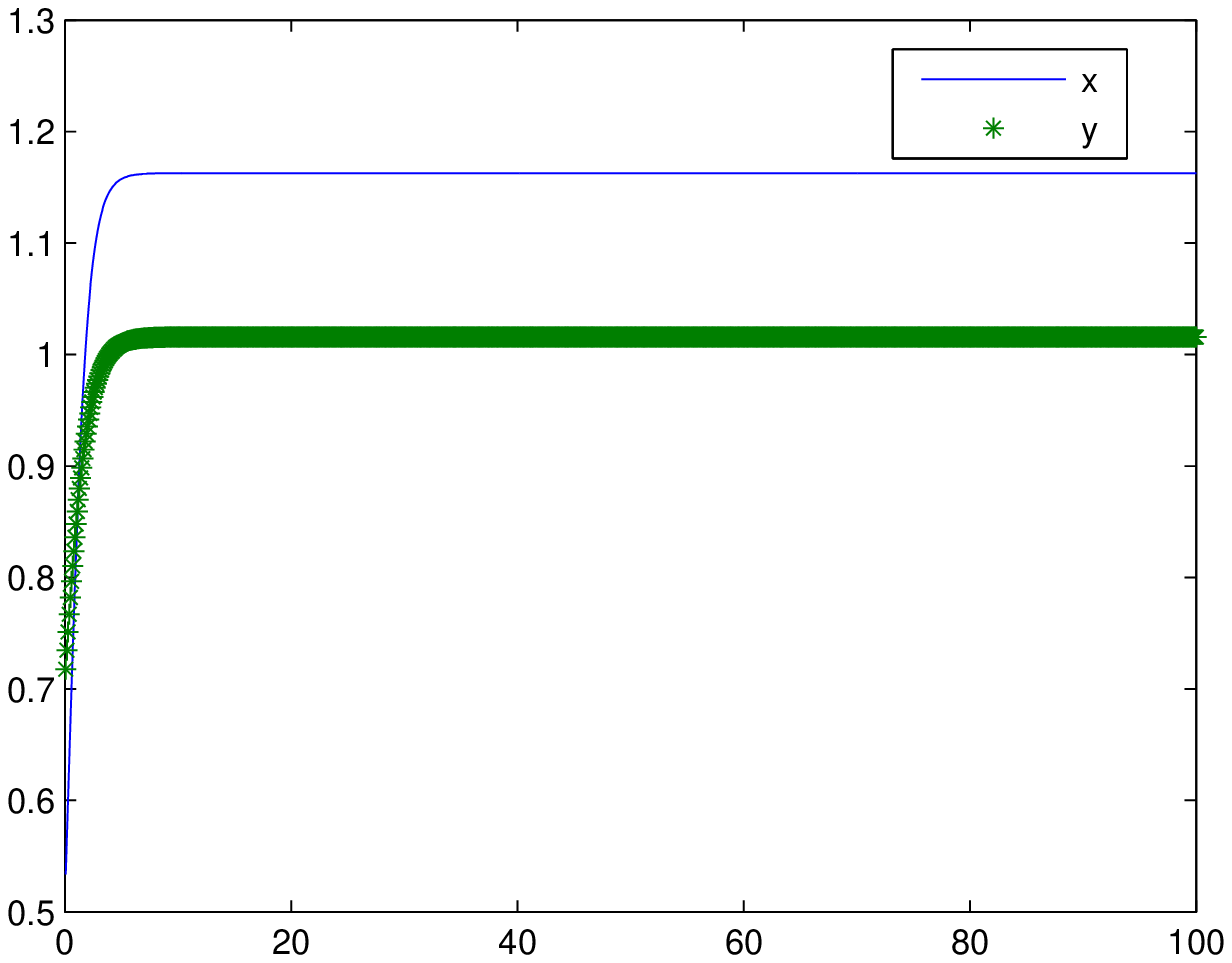}
}(b)\subfigure{
      \includegraphics[width=6.5cm]{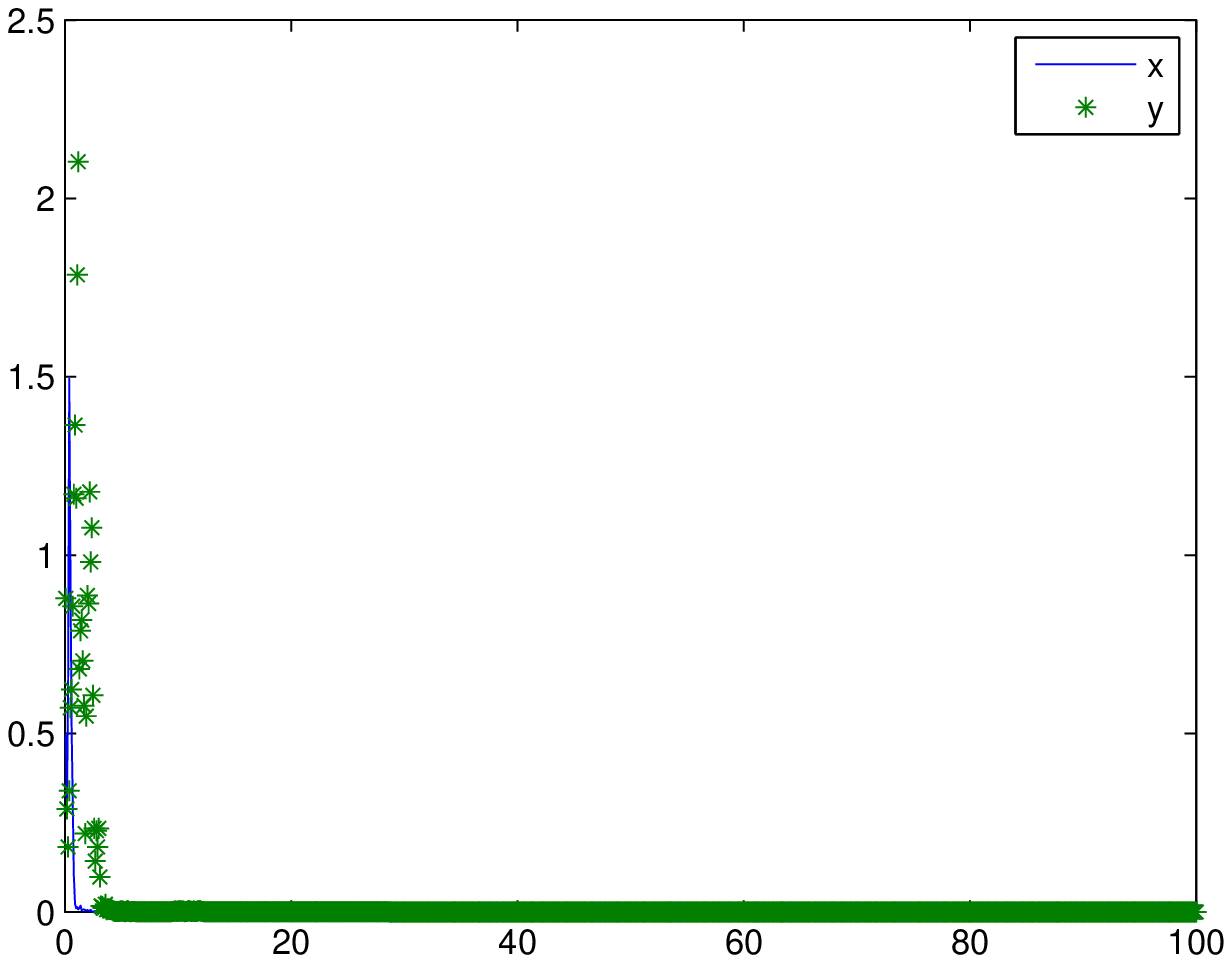}
}  (c)\subfigure{
      \includegraphics[width=6.5cm]{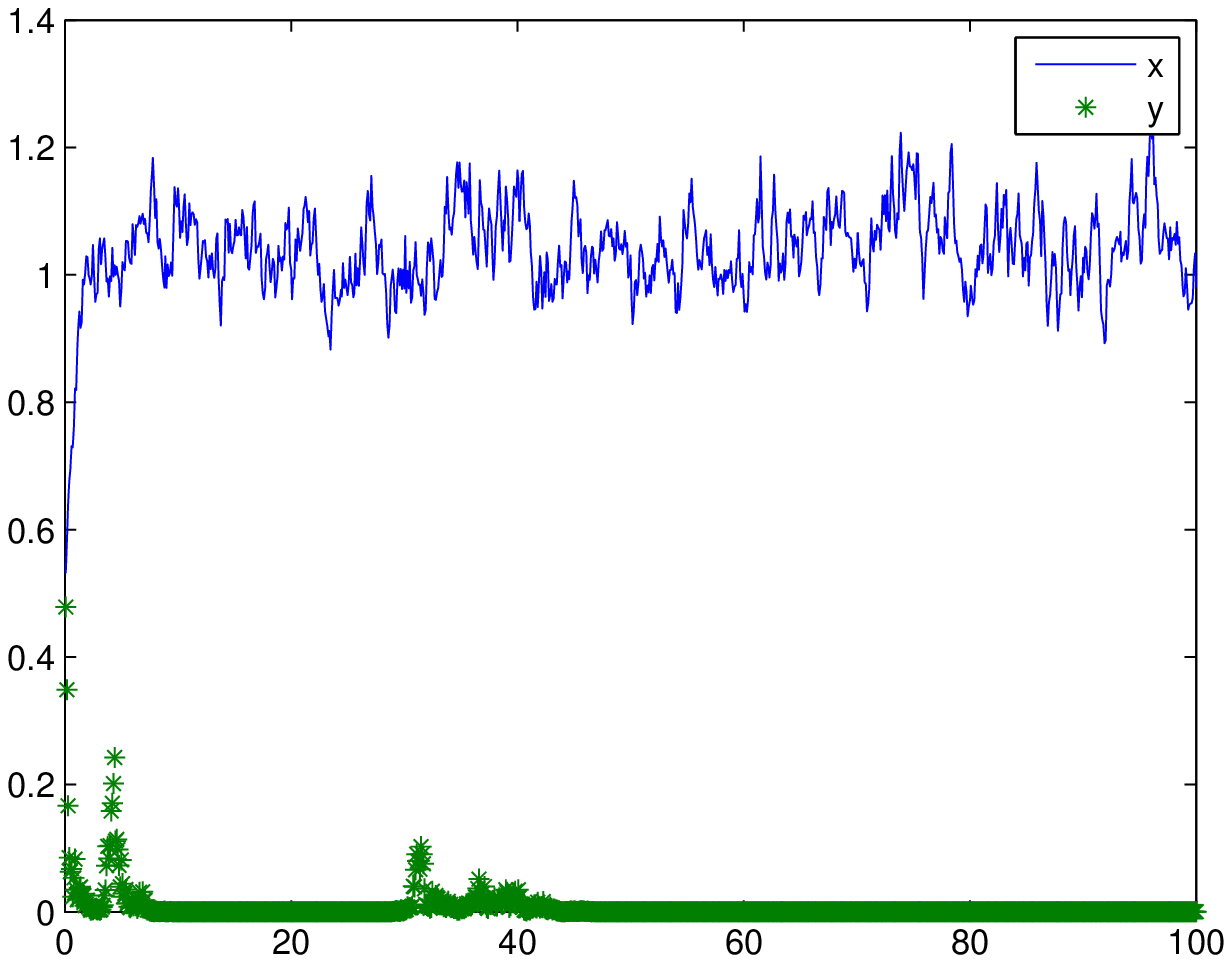}
}(d)\subfigure {
      \includegraphics[width=6.5cm]{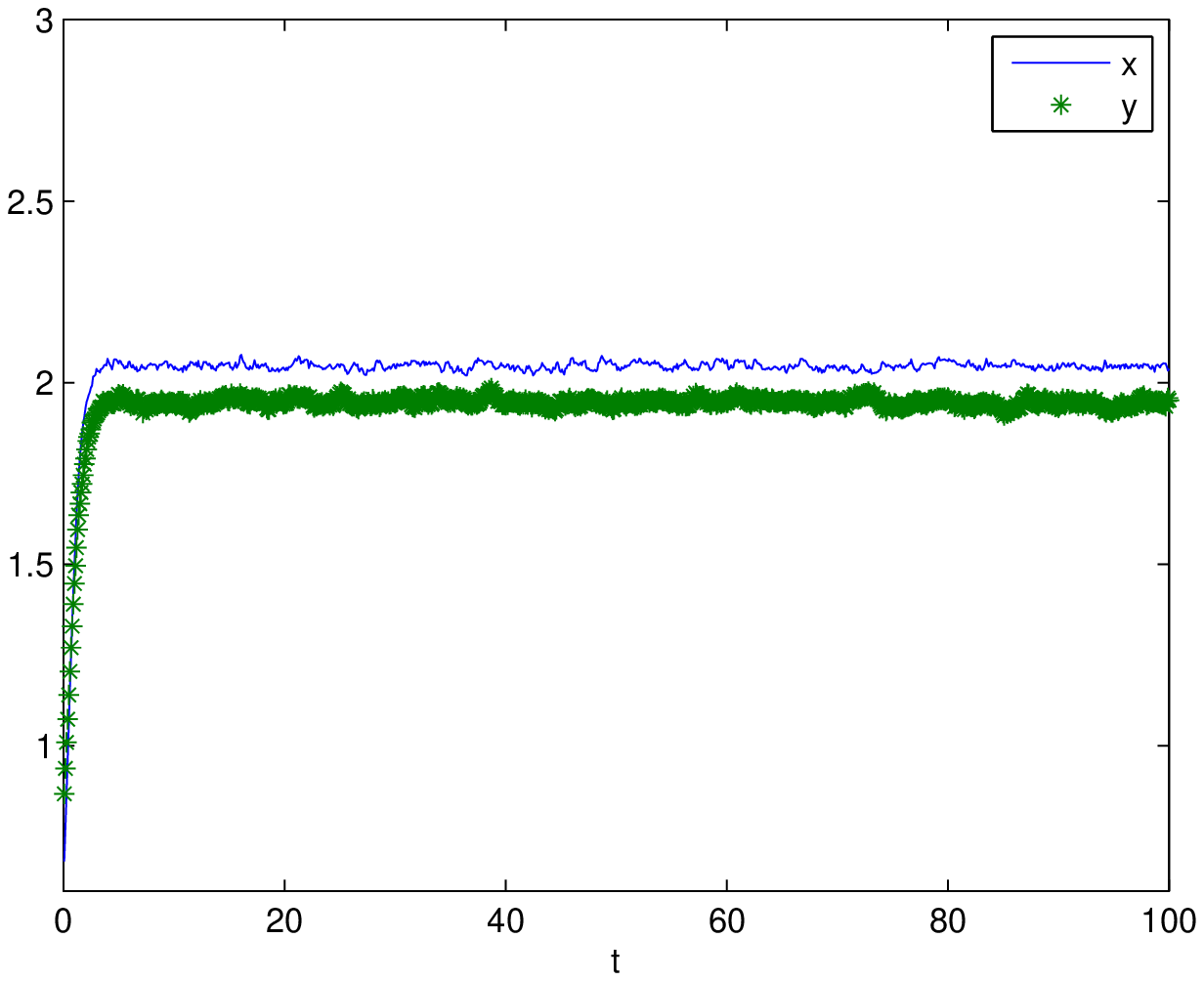}
}
 \caption{\small Solutions of \eqref{a3} for $r_1=1.2, r_2=1, \varepsilon_1=0.8, \varepsilon_2=0.7, b_1=0.7, b_2=0.9, K_1=K_2=2$, step size $\Delta$ t=0.001.
 (a) is with $\alpha_1=\alpha_2=0$; (b) is with $\alpha_1=2.2, \alpha_2=1.8 $;
 (c) is with  $\alpha_1=0.1 , \alpha_2= 1.6$; (d) is with $\alpha_1=0.01, \alpha_2=0.01.$}
\end{figure}

\end{document}